\documentclass[a4,12pt]{article}
\usepackage{graphicx}

\usepackage{tikz}
\usepackage{pgfplots}
\pgfplotsset{compat=newest, ticks=none}
\usetikzlibrary{arrows.meta}
\usetikzlibrary{calc}
\usepackage{amssymb}
\usepackage{amsmath}

\newtheorem{definition}{Definition}[section]
\newtheorem{theorem}{Theorem}[section]
\newtheorem{proposition}{Proposition}[section]
\newtheorem{corollary}{Corollary}[section]
\newtheorem{lemma}{Lemma}[section]

\date{}

\newenvironment{proof}
{\begin{rm}\par\smallskip\noindent{\bf Proof.}\quad}{\QED\end{rm}}
\def\BBox{\rule{2mm}{3mm}}
\def\QED{\hfill$\BBox$}

\def\BBox{\rule{2mm}{3mm}}
\def\QED{\hfill$\BBox$}

\title{On the Holt-Klee property for oriented matroid programming}

\author{Walter D. Morris, Jr.}
\begin{document}
\maketitle
\begin{abstract}
    The Holt-Klee theorem says that the graph of a $d$-polytope, with edges oriented by a linear function on $P$ that is not constant on any edge, admits $d$ independent monotone paths from the source to the sink.  We prove that the digraphs obtained 
    from oriented matroid programs of rank $d+1$ on $n+2$ elements, which include those from $d$-polytopes with $n$ facets, admit $d$ independent monotone paths from source to sink if $d \le 4$.  This was previously only known to hold for $d\le 3$ and $n\le 6$.  
\end{abstract}

\renewcommand{\thefootnote}{\fnsymbol{footnote}}

\section{Introduction} We will call a graph that is isomorphic to the graph formed by the vertices and edges of a $d$-dimensional convex polytope {\it d-polytopal}.  A fundamental property of $d$-polytopal graphs is that they are $d$-connected, as shown by Balinski \cite{Bal}.

The term $d$-polytopal will be applied to a digraph $K$ if there exists a $d$-dimensional polytope $P$ and an affine function $\varphi$ on $P$ that is not constant on any edge of $P$ such that the graph underlying $K$ is isomorphic to the graph of $P$ and the edges of $K$ are oriented in the direction of increase of $\varphi$.  A $d$-polytopal digraph has a unique source and sink on every subgraph corresponding to a face.  A path $(v_0,v_1,\ldots,v_k)$ from the source to the sink of a $d$-polytopal digraph will be called {\it monotone} if $(v_i,v_{i+1})$ is an arc of $K$ for $i=0,1,\ldots,k-1$.  A set of monotone paths from the source to the sink of a $d$-polytopal digraph $K$ will be called {\it independent} if the only vertices common to any 2 of the paths are the source and sink of $K$.  The following theorem comes from Holt and Klee \cite{HK}:

\begin{theorem}
 Suppose that $K$ is a $d$-polytopal digraph.  Then $K$ admits $d$ independent monotone paths from source to sink.  
  \end{theorem}  

We call the existence of $d$ independent monotone paths the Holt-Klee property.  

 Figure 1 shows acyclic orientations of 3-polytopal graphs that have a unique source and sink on every face but do not have 3 independent monotone paths from the source to the sink.  For each of the three orientations, vertices marked $A$ and $B$ cover all of the monotone paths from the source to the sink.   By the theorem of Holt and Klee, these orientations are not 3-polytopal.  The orientation of the triangular prism on the left is from Felsner, G\"artner and Tschirschnitz \cite{FGT}.  The example in the middle is an orientation of the cyclic 3-polytope with 6 vertices, from Fukuda, Moriyama and Okamoto \cite{FMO}.  The orientation of the 3-cube appeared in Stickney and Watson \cite{SW}.

\begin{figure}[h]
\setlength{\unitlength}{0.08in} % selecting unit length
\centering % used for centering Figure
\begin{picture}(65,25)
\put(5,5){\circle*{0.8}}
\put(15,5){\circle*{0.8}}\put(17,5){\makebox(0,0){A}}
\put(10,10){\circle*{0.8}}
\put(5,15){\circle*{0.8}}\put(3,15){\makebox(0,0){B}}
\put(15,15){\circle*{0.8}}
\put(10,20){\circle*{0.8}}

\put(22,14){\circle*{0.8}}
\put(38,14){\circle*{0.8}}
\put(30,6){\circle*{0.8}}\put(28,5){\makebox(0,0){A}}
\put(30,20){\circle*{0.8}}\put(32,20){\makebox(0,0){B}}
\put(26,12){\circle*{0.8}}
\put(34,12){\circle*{0.8}}

%\put(22,5){\circle*{0.8}}
%\put(32,5){\circle*{0.8}}\put(34,5){\makebox(0,0){B}}
%\put(22,15){\circle*{0.8}}\put(20,15){\makebox(0,0){A}}
%\put(32,15){\circle*{0.8}}
%\put(27,10){\circle*{0.8}}
%\put(37,10){\circle*{0.8}}
%\put(27,20){\circle*{0.8}}
%\put(37,20){\circle*{0.8}}

\put(44,5){\circle*{0.8}}
\put(54,5){\circle*{0.8}}
\put(44,15){\circle*{0.8}}\put(42,15){\makebox(0,0){A}}
\put(54,15){\circle*{0.8}}
\put(49,10){\circle*{0.8}}
\put(59,10){\circle*{0.8}}\put(61,10){\makebox(0,0){B}}
\put(49,20){\circle*{0.8}}
\put(59,20){\circle*{0.8}}

\thicklines
\put(15,5){\vector(0,1){5}}
\put(15,10){\line(0,1){5}}
\put(10,20){\vector(0,-1){5}}
\put(10,15){\line(0,-1){5}}
\put(5,15){\vector(1,0){5}}
\put(10,15){\line(1,0){5}}
\put(5,15){\vector(1,1){3}}
\put(8,18){\line(1,1){2}}
\put(5,5){\vector(1,0){5}}
\put(10,5){\line(1,0){5}}
\put(5,5){\vector(1,1){3}}
\put(8,8){\line(1,1){2}}
\put(10,10){\vector(1,-1){3}}
\put(13,7){\line(1,-1){2}}
\put(5,5){\vector(0,1){5}}
\put(5,10){\line(0,1){5}}
\put(10,20){\vector(1,-1){3}}
\put(13,17){\line(1,-1){2}}

\put(30,6){\vector(0,1){8}}
\put(30,14){\line(0,1){6}}
\put(22,14){\vector(2,-1){2}}
\put(24,13){\line(2,-1){2}}
\put(22,14){\vector(1,-1){4}}
\put(26,10){\line(1,-1){4}}
\put(22,14){\vector(4,3){4}}
\put(26,17){\line(4,3){4}}
\put(34,12){\vector(-1,0){4}}
\put(30,12){\line(-1,0){4}}
\put(34,12){\vector(2,1){2}}
\put(36,13){\line(2,1){2}}
\put(26,12){\vector(1,2){2}}
\put(27,14){\line(1,2){3}}
\put(30,6){\vector(-2,3){2}}
\put(28,9){\line(-2,3){2}}
\put(30,6){\vector(2,3){2}}
\put(32,9){\line(2,3){2}}
\put(30,6){\vector(1,1){4}}
\put(34,10){\line(1,1){4}}
\put(34,12){\vector(-1,2){2}}
\put(33,14){\line(-1,2){3}}
\put(30,20){\vector(4,-3){4}}
\put(34,17){\line(4,-3){4}}

%\put(22,5){\vector(1,0){10}}
%\put(22,5){\vector(0,1){10}}
%\put(22,5){\vector(1,1){5}}
%\put(32,5){\vector(0,1){10}}
%\put(22,15){\vector(1,0){10}}
%\put(37,20){\vector(-1,-1){5}}
%\put(27,20){\vector(1,0){10}}
%\put(37,20){\vector(0,-1){10}}
%\put(22,15){\vector(1,1){5}}
%\put(27,20){\vector(0,-1){10}}
%\put(27,10){\vector(1,0){10}}
%\put(37,10){\vector(-1,-1){5}}

\put(44,5){\vector(1,0){5}}
\put(49,5){\line(1,0){5}}
\put(44,5){\vector(0,1){5}}
\put(44,10){\line(0,1){5}}
\put(44,5){\vector(1,1){3}}
\put(47,8){\line(1,1){2}}
\put(44,15){\vector(1,0){5}}
\put(49,15){\line(1,0){5}}
\put(49,10){\vector(1,0){5}}
\put(54,10){\line(1,0){5}}
\put(54,5){\vector(1,1){3}}
\put(57,8){\line(1,1){2}}
\put(44,15){\vector(1,1){3}}
\put(47,18){\line(1,1){2}}
\put(49,20){\vector(0,-1){5}}
\put(49,15){\line(0,-1){5}}
\put(54,15){\vector(0,-1){5}}
\put(54,10){\vector(0,-1){5}}
\put(49,20){\vector(1,0){5}}
\put(54,20){\line(1,0){5}}
\put(54,15){\vector(1,1){3}}
\put(57,18){\line(1,1){2}}
\put(59,10){\vector(0,1){5}}
\put(59,15){\line(0,1){5}}
\end{picture}
\caption{Only two independent monotone paths from source to sink} % title of the Figure
\label{fig:lnlblock} % label to refer figure in text
\end{figure}
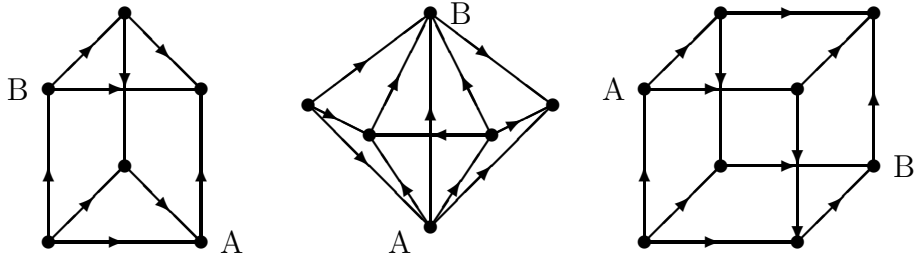

Combinatorial types of polytopes are difficult to enumerate because non-obvious restrictions quickly arise as the dimension and number of faces increase.  Oriented matroids provide a way to generate combinatorial spheres that satisfy a short list of necessary conditions for polytopes.  There are two dual spheres that may be derived from an oriented matroid.  The Las Vergnas sphere has vertices corresponding to the elements of the oriented matroid ground set.  The Edmonds - Mandel sphere has facets corresponding to elements of the ground set.  The papers by Bremner et. al. (\cite{BremSchew}, \cite{BDHS}) use the term  matroid polytope for the Las Vergnas sphere of an oriented matroid.  They studied ridge paths of the Las Vergnas sphere, which correspond to paths in the graph of the Edmonds - Mandel sphere.  Their upper bounds on the diameter of polytopal graphs for fixed small dimension and number of facets were obtained by showing that the bound holds for graphs of Edmonds - Mandel spheres.  It is known that Edmonds - Mandel spheres are shellable, but shellability of Las Vergnas spheres 
is an open problem.

An {\it oriented matroid program} is an oriented matroid ${\mathcal M}$ of rank $d+1$ on a ground set of size $n + 2$, where $n$ elements of the ground set correspond to facets of the $(d-1)$-dimensional Edmonds - Mandel sphere. The remaining two elements, called the objective and right hand side elements, determine the orientation of the graph of the Edmonds - Mandel sphere.  If ${\mathcal M}$ is realizable, these 
define an orientation of the graph of a polytope that is $d$-polytopal.  Oriented matroid programming has had a central place in the evolution of oriented matroid theory, starting with \cite{Bland}, \cite{FL}, and continuing with \cite{EF}, \cite{Todd} and \cite{KO}.

\begin{figure}[h]    
\setlength{\unitlength}{0.10in} % selecting unit length
\centering % used for centering Figure
\begin{picture}(45,25)

\put(14,5){\circle*{0.8}}
\put(24,5){\circle*{0.8}}
\put(14,15){\circle*{0.8}}\put(13,4){\makebox(0,0){v}}
\put(24,15){\circle*{0.8}}
\put(19,10){\circle*{0.8}}
\put(29,10){\circle*{0.8}}\put(30,21){\makebox(0,0){w}}
\put(19,20){\circle*{0.8}}
\put(29,20){\circle*{0.8}}

\thicklines
\put(14,5){\vector(1,0){5}}
\put(14,5){\vector(0,1){5}}
\put(14,5){\vector(1,1){3}}
\put(24,15){\vector(-1,0){5}}
\put(19,10){\vector(1,0){5}}
\put(29,10){\vector(-1,-1){3}}
\put(14,15){\vector(1,1){3}}
\put(19,20){\vector(0,-1){5}}
\put(24,5){\vector(0,1){5}}
\put(19,20){\vector(1,0){5}}
\put(24,15){\vector(1,1){3}}
\put(29,10){\vector(0,1){5}}

\put(19,5){\line(1,0){5}}
\put(14,10){\line(0,1){5}}
\put(17,8){\line(1,1){2}}
\put(19,15){\line(0,-1){5}}
\put(24,10){\line(1,0){5}}
\put(26,7){\line(-1,-1){2}}
\put(17,18){\line(1,1){2}}
\put(19,15){\line(-1,0){5}}
\put(24,10){\line(0,1){5}}
\put(24,20){\line(1,0){5}}
\put(27,18){\line(1,1){2}}
\put(29,15){\line(0,1){5}}
\end{picture}

\caption{Oriented matroid program with a monotone cycle} % title of the Figure
\label{fig:lnlblock} % label to refer figure in text
\end{figure}
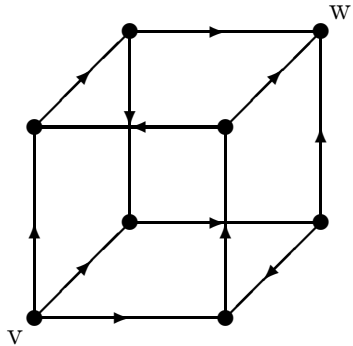

Figure 2 shows that the extension of polytopal digraphs to digraphs of oriented matroid programs is proper already for $d=3$ and $n=6$, coming from an oriented matroid of rank 4 on 8 elements (See \cite{EF}, \cite{EM}). The digraph depicted cannot be 3-polytopal, due to the monotone cycle
containing all vertices other than the source and sink.

%Much of the early research on oriented programming aimed at finding 
%a simplex algorithm for oriented matroid programming that avoided directed cycles and maintained feasibility once attained.  These goals were attained by M. J. Todd \cite{Todd}.  The path followed by Todd's algorithm once feasibility is attained is a generalization of a {\it coherent path}. 

By Menger's theorem \cite{Menger}, a $d$-polytopal digraph admits $d$ independent monotone paths from source to sink if and only if there is no set of $d-1$ vertices, not containing the source and the sink, that covers all of the monotone source to sink paths.  Given a $d$-dimensional polytope $P$ and a set $S$ of $d-1$ vertices not containing the source and sink, the argument of Holt and Klee \cite{HK} introduces a hyperplane $H$, containing the source and sink, for which one of the open half-spaces $H^-$ determined by $H$ contains points of $P$ but no points of $S$.  Then a monotone source to sink path is found in the union of $H^-$ with the set consisting of the source and sink. The extension of this proof to oriented matroid programming would require the extension of the oriented matroid program by an element corresponding to the  hyperplane $H$. There are several known obstructions to extending a non-realizable rank 4 oriented matroid. It is therefore reasonable that  
the seminal paper \cite{FMO} on the Holt-Klee property for oriented matroids showed by computer enumeration that every oriented matroid program of rank 4 on 8 elements has the Holt - Klee property, that is, the orientation of the graph of the  Edmonds - Mandel sphere has 3 independent monotone paths from the source to the sink. 
We adapt the argument of \cite{HK} and \cite{GMR} to the setting of oriented matroids to 
prove that an oriented matroid program ${\mathcal M}$ has the Holt-Klee property if an extension problem in the oriented matroid ${\mathcal M}/\{f,g\}$, of rank 2 less than that of ${\mathcal M}$, has a solution.  This will, in particular, prove that ${\mathcal M}$ has the Holt-Klee property if $d=3$ or $d=4$, with $n$ arbitrarily large.

The $d=3$ case has particular interest because \cite{HK} used the Holt-Klee property to prove the bound $\Delta_{sm}(3,n) \le \lfloor \frac{2}{3}n\rfloor-1$ on the strict monotone diameter in dimension 3.  Our results imply that this bound also holds for oriented matroid programs. See \cite{GM} for a recent application of oriented matroid programming to determination of polytope diameters. In dimension 3, the Holt-Klee condition is also known to be sufficient for an acyclic orientation of the graph of a polytope with a unique source and sink on every face to be $d$-polytopal. (See \cite{MK},\cite{LM}).

The paper \cite{FMO} showed that oriented matroid programming is not the only way to produce digraphs from oriented matroids that are $d$-polytopal in the realizable case.  From an oriented matroid of rank $d+1$ on $n + (d-1)$ elements, where $n$ of the elements correspond to the facets of the Edmonds - Mandel sphere and $d-1$ define a line shelling of the sphere, one obtains an orientation of the graph of the {\it Las Vergnas} sphere.  Their enumeration of rank 4 oriented matroids on 8 elements produced the digraph in the middle of Figure 1, which does not have 3 independent monotone paths from the source to the sink.  It follows that digraphs produced this way from oriented matroids do not necessarily satisfy the Holt-Klee condition, even though those that are produced from realizable oriented matroids do.

\section{Oriented Matroids }  We will give some basic definitions, following the notation of the standard reference \cite{OMbook}.  

Let $E$ be a finite set.  A {\it sign vector} on $E$ is a vector indexed by $E$ with entries in 
$\{+,-,0\}$. The {\it support} of a sign vector $X$ is $\underline{X}:=\{e \in E: X_e \neq 0\}$.  For a sign vector $X$, we define $X_+=\{e\in E: X_e=+\}$ and $X_-=\{e\in E: X_e=-\}$.   The {\it negative} of a sign vector 
$X$ is the sign vector $-X$ obtained from $X$ by replacing $+$ entries with $-$ and vice versa.  The {\it composition} $X \circ Y$ of sign vectors $X$ and $Y$ is defined by $(X \circ Y)_e = X_e$ if $X_e \neq 0$ and $(X \circ Y)_e = Y_e$ otherwise.  For sign vectors $X,Y$, define $S(X,Y) = \{e \in E: X_e = -Y_e \neq 0\}$.  Two sign vectors $X$ and $Y$ are said to be {\it orthogonal} if $S(X,Y)$ and $S(X,-Y)$ are both empty or both nonempty.  

An {\it oriented matroid} ${\mathcal M}$ on $E$ consists of a collection ${\mathcal Y}$ of sign vectors, called {\it covectors}, on $E$, satisfying the following properties:  

\begin{description}
\item [(Y0)] $0 \in {\mathcal Y}$,
\item [(Y1)] If $X \in {\mathcal Y}$ then $-X \in {\mathcal Y}$,
\item [(Y2)] If $X,Y \in {\mathcal Y}$ then $X \circ Y \in {\mathcal Y}$,
\item [(Y3)] If $X,Y \in {\mathcal Y}$, and $e \in S(X,Y)$, then there exists $Z \in {\mathcal Y}$ so that $Z_e = 0$ and $Z_f = (X \circ Y)_f = (Y \circ X)_f$ for all $f \in E\backslash S(X,Y)$.
\end{description}

A covector $Z$ obtained as in $(Y3)$ is said to be obtained by {\it eliminating} $e$ between $X$ and $Y$.  
The nonzero covectors in ${\mathcal Y}$ with minimal support are called the {\it cocircuits} of ${\mathcal M}$.  The set of cocircuits of ${\mathcal M}$ is denoted ${\mathcal C}^*({\mathcal M})$.  We will call a covector $Y$ {\it nonnegative} if $Y_e \ge 0$ for all $e \in E$.  

A sign vector on $E$ is called a {\it vector} of ${\mathcal M}$ if it is orthogonal to every covector of ${\mathcal M}$.  The nonzero vectors of ${\mathcal M}$ with minimal support are called {\it circuits} of ${\mathcal M}$.  We denote by ${\mathcal C}({\mathcal M})$ the set of circuits of ${\mathcal M}$.  The vectors of ${\mathcal M}$ are the covectors of an oriented matroid ${\mathcal M}^*$ on $E$, called the {\it dual} of ${\mathcal M}$.  

A subset of $E$ is called {\it independent} if it does not contain the support of a circuit.  The {\it rank} of ${\mathcal M}$ is the size of a largest independent subset of $E$. A maximal independent subset of ${\mathcal M}$ is a {\it base} of ${\mathcal M}$.  A base of ${\mathcal M}^*$ is called a {\it cobase} of ${\mathcal M}$.
 
If $B$ is a base of ${\mathcal M}$ and $p \in E\backslash B$, then we denote by $X(B,p)$ the unique circuit of ${\mathcal M}$ with support in $B \cup \{p\}$ and with coordinate $p$ equal to $+$.  The complement ${\overline B} = E\backslash B$ is a cobase of ${\mathcal M}$.
For every $q \in B$, there is a unique cocircuit $Y(\overline{B},q)$ with support in $\overline B \cup q$ and with $Y(\overline{B},q)_q =+$.

If $p \in E$ and $A \subseteq E\backslash p,$
we say $p \in conv_{\mathcal M}(A)$ if there is a circuit $X$ of ${\mathcal M}$ such that  $X_+ \subseteq A$ and $X_-=\{p\}$.  We also say that $-p \in conv_{\mathcal M}(A)$ if there is a positive circuit $X$ of ${\mathcal M}$ such that $p \in X_+\subseteq A \cup \{p\}$ and $X_-=\emptyset.$  

If $T$ is a matrix with real entries, then the collection of sign vectors obtained from the vectors in the row space of $T$  by forgetting magnitudes is the set of vectors of an oriented matroid, called the oriented matroid {\it realized} by $T$. The covectors of this oriented matroid correspond analogously to vectors in the null space of $T$.

If $X$ is a sign vector on $E$ and $e \in E$, we denote by $X \backslash e$ the subvector of $X$ with component indexed by $e$ removed.  
If ${\mathcal M}$ is an oriented matroid on $E$ and $e \in E$, the oriented matroid ${\mathcal M}\backslash e$ on $E\backslash e$ has circuit set equal to $\{X \backslash e: X \in {\mathcal C}({\mathcal M}), X_e = 0\}$.  We say that ${\mathcal M}\backslash e$ is obtained from ${\mathcal M}$ by deleting $e$.  If $A \subseteq E$ then ${\mathcal M}(A)$, the restriction of ${\mathcal M}$ to $A$, is the oriented matroid obtained from ${\mathcal M}$ by deleting $E\backslash A$.  The oriented matroid ${\mathcal M}/e$ on $E\backslash e$ has circuit set the members of $\{X\backslash e: X \in {\mathcal C}({\mathcal M})\}$ with minimal support. ${\mathcal M}/e$ is said to be obtained from ${\mathcal M}$ by contracting $e$. The relationship  $({\mathcal M}\backslash e)^* = {\mathcal M}^*/e$ shows how duality and minors are related.

If ${\mathcal M}$ is an oriented matroid on a set $E$ and $\hat{e} \notin E$, then an oriented matroid $\hat{\mathcal M}$ on $E \cup \hat{e}$  is said to be an extension of ${\mathcal M}$ to $E \cup \hat{e}$ if $\hat{\mathcal M} \backslash \hat{e} = {\mathcal M}$.  We will only consider extensions for which the rank of $\hat{\mathcal M}$ equals the rank of 
${\mathcal M}$.  An oriented matroid $\hat{\mathcal M}$ such that $(\hat{\mathcal M})^*$ is an extension of ${\mathcal M}^*$ will be called a {\it dual extension} of ${\mathcal M}$.

Suppose ${\mathcal M}$ is an oriented matroid on a set $E$, and $\sigma$ is a function from ${\mathcal C}^*({\mathcal M})$, the set of cocircuits of ${\mathcal M}$, to $\{+,-,0\}$, satisfying $\sigma(-Y) = -\sigma(Y)$ for every cocircuit $Y$ of ${\mathcal M}$. For every cocircuit $Y \in {\mathcal C}^*({\mathcal M})$, define $\hat{Y}$ on $E \cup e$ so that $\hat{Y}_f = Y_f$ for all $f \in E$ and $\hat{Y}_e = \sigma(Y)$. Theorem 7.1.8 of \cite{OMbook}, due to Las Vergnas, says that $\{\hat{Y}:Y \in {\mathcal C}^*({\mathcal M})\}$ is a subset of the set of cocircuits of an extension of ${\mathcal M}$ if and only if $\sigma$ defines a single element extension of every contraction of ${\mathcal M}$ of rank 2.

\section{Oriented matroid programming}

Oriented matroid programming provides a combinatorial abstraction of linear programs of the form

$$  {\begin{array}{*{20}c}
 \text{maximize } f = c^Ty + \beta \text{ subject to}  \\
    Ay = b \\
    y \ge 0.
      \end{array} }. $$

Here $A$ is assumed to be an $m \times n$ matrix of rank $m$.  We view linear programming as the search for vectors of certain sign patterns in the null space and row space of the matrix 

$$
T = \left[
\begin{array}{c|c|c}
A & 0& -b \\ \hline
 -c^T& 1&-\beta
\end{array}\right].
$$

With this motivating linear program in mind, we assume that ${\mathcal M}$ is a rank $r=n-m+1$ oriented matroid on an $(n+2)$-element  set $E = [n] \cup \{f,g\}$.  The element $f$ is called the objective element and the element $g$ is called the right hand side element.

When  ${\mathcal M}$ is realized by $T$, the covectors of ${\mathcal M}$ are sign vectors for vectors in the null space of $T$.   The elements $f$ and $g$ correspond to the last two columns of $T$.  The covectors of ${\mathcal M}\backslash f$ are obtained by deleting component $f$ from the covectors of ${\mathcal M}$.  Covectors $Y$ of ${\mathcal M}\backslash f$ with $Y_g = +$ correspond to solutions of $Ay=b$ in the realizable case.  

%We will assume that $f$ is not a coloop, so that the rank of ${\mathcal M}\backslash f$ equals the rank of ${\mathcal M}$.  A covector $Y$ of the oriented matroid ${\mathcal M}\backslash f$ is called {\it feasible} if $Y_g = +$ and $Y_e \ge 0$ for all $e \in [n]$.  A cobase $\overline{B}$ of ${\mathcal M}$ will be called $\it feasible$ if $f \in \overline{B}, g \notin \overline{B}$, and $Y(\overline{B},g)\backslash f$ is a feasible cocircuit of ${\mathcal M}\backslash f$.   A cocircuit $Y$ of ${\mathcal M}$ is called feasible if $Y = Y(\overline{B},g)$ for some feasible cobase $\overline{B}$.    

\subsection{The underlying undirected graph}

 We will assume that every nonnegative cocircuit of ${\mathcal M}\backslash f$ contains $g$ in its support.   In the terminology of \cite{OMbook}, the feasible region of the oriented matroid program defined by ${\mathcal M}$ contains no covectors at infinity.  This implies the oriented matroid program is {\it bounded}.  We will also assume that for every $e \in [n]$ there is a nonnegative cocircuit $Y$ of ${\mathcal M}\backslash f$ with $Y_e = +$.  In the terminology of \cite{OMbook}, the oriented matroid ${\mathcal M}\backslash f$ is {\it acyclic}.

Consider the set ${\mathcal F} := \{[n]\backslash \underline{Y}: Y \mbox{ is a nonnegative covector of }{\mathcal M}\backslash f\}$.  The collection ${\mathcal F}$ can be partially ordered, with 
$F_1 \le F_2$ if $F_1 \subseteq F_2$.  The elements of ${\mathcal F}$ are the {\it faces} of the oriented matroid polytope given by ${\mathcal M}\backslash f$, 
and the partial order is called the Las Vergnas face lattice.  Its order dual is called the Edmonds-Mandel face lattice.  
\begin{theorem}(Theorem 4.3.5 of \cite{OMbook}) Let $r$ be the rank of ${\mathcal M}$ and ${\mathcal M}\backslash f$.   \begin{enumerate}\item The Edmonds-Mandel lattice is isomorphic to the face lattice of a shellable regular cell decomposition of the $(r-2)$-sphere. 
\item The Las Vergnas lattice is isomorphic to the face lattice of a PL regular cell decomposition of the $(r-2)$-sphere.  
\end{enumerate} 
\end{theorem}

We are interested in the graph $G$ of the Edmonds-Mandel sphere.  We will use the words {\it facet} and {\it node} for an element of ${\mathcal F}$ of rank $r-1$, depending on whether it is viewed as a facet of the Las Vergnas sphere or a vertex of the Edmonds-Mandel sphere.  
Two distinct nodes $v$ and $v'$ of $G$ are adjacent if $v \cap v'$ has rank $r-2$.     
The graph is $(r-1)$-connected, by Theorem 4.4.9 of \cite{OMbook}.  

Because of the assumption that $Y_g = +$ for every nonnegative covector of ${\mathcal M}\backslash f$, it follows that the poset of nonnegative covectors of ${\mathcal M}\backslash \{f,g\}$ ordered as before, is the same as ${\mathcal F}$.  Thus the elements $f$ and $g$ are not required to define the undirected graph $G$.  

%If we make the additional assumption that for every $e \in [n]$ there is a positive covector $Y$ of ${\mathcal M} \backslash \{f,g\}$ for which $e$ is the only component equal to $0$, then the Las Vergnas lattice is the face lattice of a {\it matroid polytope}.  

\subsection {The orientation of the graph}
The digraph $K_f$ of the oriented matroid program is obtained by orienting the edges of $G$. We will assume that every circuit of ${\mathcal M}$ containing $f$ in its support has support of size $r+1$.  This will ensure that every edge of $G$ receives an orientation.  

Every node of $K_f$ is the zero set of a cocircuit  $Y$ of ${\mathcal M}$ that is nonnegative on $E\backslash f$, and has $f,g \in \underline{Y}$.  Suppose $v$ and $v'$ corresponding to $Y$ and $Y'$ are adjacent nodes of $K_f$.  Let the cocircuit $Y''$ be obtained by eliminating $g$ between $Y$ and $-Y'$.  Then removing the zero entry in position $g$ from $Y''$ yields a cocircuit of ${\mathcal M}/g$.  Note that $Y''_e = +$ for all $e$ in $v'\backslash v$ and $Y''_e=-$ for all $e$ in $v\backslash v'$, and $Y''_e=0$ for all $e$ in $v \cap v'$. %Then the restriction of the vector $-Y''$ to the entries indexed by $v$ %is a positive covector of ${\mathcal M}(v)$ and the restriction of $Y''$ to %the entries indexed by $v'$ is a positive covector of ${\mathcal M}(v')$.  
Because of the nondegeneracy assumption on $f$ we have $Y''_f \neq 0$. If $Y''_f = Y''_e$ for $e \in v\backslash v'$, then the arc containing $v$ and $v'$ is directed from $v$ to $v'$.

If $B$ is a base of ${\mathcal M}$ contained in the set $([n] \backslash \underline{Y''}) \cup \{e,g\}$ for some $e \in v\backslash v'$, then the  supports of $X=C(B,f)$ and $Y''$ intersect on $\{e,f\}$ Thus the arc leaves $v$ if entry $e$ of $X$ is negative. 

The fundamental theorem of oriented matroid programming (10.1.13 in \cite{OMbook}) implies that the digraph $K_f$ has a unique sink.  By considering minors, it implies that the restriction of $K_f$ to any face of ${\mathcal M}\backslash f$ has a unique sink.   
   
We summarize the above discussion:
The undirected graph $G$ of the oriented matroid program is defined by the cocircuits of the minor ${\mathcal M}\backslash f$.
The edges of the graph are given orientations by the minor ${\mathcal M} / g$.

\section{Building a Path from Source to Sink}

The path that we will build is in the nondegenerate case a generalization to oriented matroids of the Gass-Saaty shadow-vertex algorithm.  For cases when the Edmonds-Mandel sphere is not simple, it is most succinctly described in section 4 of \cite{S}.   This description is in terms of the cells of a subdivision of an oriented matroid. We repeat Definition 4.5 of \cite{S}:

\begin{definition} \label{defsubdivision}A collection $S$ of subsets of $E$ (called {\it cells}) is called a {\it subdivision} of a rank $r-1$ oriented matroid ${\mathcal N}$ on $E$ if the following hold:  
\begin{enumerate}
\item Every cell $\sigma \in S$ has rank $r-1$. (The ``$r$" in \cite{S} is replaced by $r-1$ here, the rank of ${\mathcal M}\backslash f/g$.)
\item For every one-element extension ${\mathcal N} \cup f$ of ${\mathcal N}$ 
and every $\sigma_1,\sigma_2 \in S$, 
$$f \in conv_{{\mathcal N} \cup f}(\sigma_1) \cap conv_{{\mathcal N} \cup f}(\sigma_2) \Rightarrow f \in conv_{{\mathcal N} \cup f}(\sigma_1 \cap \sigma_2)$$
\item If $\sigma_1,\sigma_2 \in S$, then $\sigma_1 \cap \sigma_2$ is a common face of the two restrictions ${\mathcal N}(\sigma_1)$ and ${\mathcal N}(\sigma_2)$.
\item If $\sigma \in S$, then each facet of ${\mathcal N}(\sigma)$ is either contained in a facet of ${\mathcal N}$ or contained in precisely two cells of $S$.  
\end{enumerate}
\end{definition}

A subdivision is called a {\it lifting} subdivision if there is an oriented matroid $\widehat{\mathcal N}$ on a set $E \cup g$ so that ${\mathcal N}= \widehat{\mathcal N}/g$ and the subsets in $S$ are the complements of the supports of nonnegative cocircuits of $\widehat{\mathcal N}$.   

We will be concerned with the lifting subdivision of ${\mathcal M}\backslash f/g$ given by the facets of the Las Vergnas sphere of ${\mathcal M}\backslash f$.  Recall that the vertices of the Edmonds-Mandel sphere are the facets of the Las Vergnas sphere.  

The oriented matroid ${\mathcal M}/g$ is an extension of ${\mathcal M}\backslash f/g$ by an element $f$ in general position.  If $v$ is the unique sink of the digraph $K_f$, then $-f \in conv_{{\mathcal M}/g}(v)$.  By considering the oriented matroid program obtained from ${\mathcal M}$ by reversing the sign of $f$ in every vector and covector of ${\mathcal M}$, one sees that $K_f$ has a unique source.  If $v'$ is the unique source of $K_f$, then $f \in conv_{{\mathcal M}/g}(v')$.  
With the interpretation of vertices of $K_f$ as facets of the Las Vergnas sphere, or as cells of a subdivision, we say that $-f$ is covered by $v$ and $f$ is covered by $v'$.  

Now we will extract from Lemma 4.6 of \cite{S} what we need. Define $\widehat{{\mathcal M}}$ to be an extension of ${\mathcal M}/g$ by an element $h$ in general position.  Then $\widehat{{\mathcal M}}$ is a rank $r-1$ oriented matroid on $[n] \cup \{f,h\}$ satisfying $\widehat{{\mathcal M}}\backslash h = {\mathcal M}/g.$ Because $h$ is in general position, every circuit of $\widehat{{\mathcal M}}$ containing $h$ has $r$ elements.  

We first note that one can define an alternate orientation $K_h$ of the graph of the Edmonds-Mandel sphere using $\widehat{{\mathcal M}}\backslash f$.  Suppose $v$ and $v'$ corresponding to $Y$ and $Y'$ are adjacent nodes of $K_f$.  Recall the cocircuit $Y''$ of ${\mathcal M}/g$ that was used in defining the orientation of the arc between $v$ and $v'$.  If we remove entry $Y''_f$ from $Y''$, we get a cocircuit of ${\mathcal M}\backslash f/g$ which agrees with $Y''$ on $[n]$.  This cocircuit corresponds to a cocircuit $\widehat{Y}''$ of $\widehat{{\mathcal M}}\backslash f$ with an extra entry indexed by $h$.  If $\widehat{Y}''_h = \widehat{Y}''_e$ for $e \in v\backslash v'$, then the arc containing $v$ and $v'$ is directed from $v$ to $v'$ in $K_h$.

\begin{lemma}
The orientation $K_h$ has a unique sink.
\end{lemma}
\begin{proof}
This follows from part (i) of Lemma 4.6 of \cite{S}.  $\widehat{{\mathcal M}}\backslash f$ is an extension of ${\mathcal M}\backslash f/g$ by the element $h$ in general position, so $h$ is covered by a unique facet of the Las Vergnas sphere of $\widehat{{\mathcal M}}\backslash f$.  By considering the oriented matroid obtained from $\widehat{{\mathcal M}}\backslash f$ by reversing the sign of $h$, it follows that $-h$ is also covered by a unique facet.  
\end{proof}

Each of the oriented matroids $\widehat{{\mathcal M}}\backslash f$ and $\widehat{{\mathcal M}}\backslash h = {\mathcal M}/g$ is an extension of ${\mathcal M}\backslash f/g$ by an element in general position.  The reason for using Lemma 4.6 of \cite{S} (rather than the fundamental theorem of oriented matroid programming) to establish that $K_h$ has a unique sink is that there is not necessarily an oriented matroid program ${\mathcal M}'$ on $[n] \cup \{h,g\}$ for which ${\mathcal M}'/g=\widehat{{\mathcal M}}\backslash f$.  The extension of ${\mathcal M}\backslash f/g$ by $h$ and the dual extension by $g$ may be incompatible.  

The proof in \cite{S} of Lemma 4.6 goes on to construct a digraph $G_{[p_1,p_2]}$, which we will call $K_{[f,h]}$:

\begin{definition}\label{Gfhdef}
Let ${\mathcal N}'={\mathcal N}\cup \{f,h\}$ be a two-element extension of ${\mathcal N}$, where each of ${\mathcal N}'\backslash f$ and ${\mathcal N}'\backslash h$ is an extension of ${\mathcal N}$ by an element in general position.  Form the digraph $K_{[f,h]}$ as follows. \begin{itemize}
    \item A cell $\sigma \in S$ is a vertex of the graph if and only if there is a vector $X$ of ${\mathcal N}'$ with $X_+ = \sigma$ and $X_- = \{f,h\}$.
    \item Let $\tau$ be an $(r-2)$-face of a cell of $S$ for which there is a vector $X$ of ${\mathcal N}'$ with $X_+ = \tau$ and $X_- = \{f,h\}.$  Then there are exactly two cells $\sigma^+$ and $\sigma^-$ of $S$ containing $\tau$.  Let $Y$ be the cocircuit of ${\mathcal N}'$ vanishing on $\tau.$  Then $Y_{f}=-Y_{h}$. 
    Introduce an arc oriented from the vertex $\sigma^+$ to the vertex $\sigma^-$ if $f \cup (\sigma^+ \backslash \tau) \subseteq Y^+$ and  $h \cup (\sigma^- \backslash \tau) \subseteq Y^-$.
\end{itemize}
\end{definition}  

In our application, let ${\mathcal N}$ be ${\mathcal M}\backslash f/g$, let ${\mathcal N}'$ be $\widehat{{\mathcal M}}$, and let $S$ be the set of facets of the Las Vergnas sphere of ${\mathcal M}\backslash f$.  Then $K_{[f,h]}$ can be seen to be a subdigraph of $K_f$.  

\begin{lemma}  (From the proof of Lemma 4.6 of \cite{S})
The component of $K_{[f,h]}$ containing the source of $K_f$ is a monotone path of $K_f$ from the source of $K_f$ to the source of $K_h$.  
\end{lemma}

We can adapt Definition \ref{Gfhdef} to construct two subdigraphs of $K_f$ which we call $K_{[-h,-f]}$ and $K_{[f,-h]}.$  A node $v$ of $K_f$ is a node of $K_{[-h,-f]}$ if there is a nonnegative vector $X$ of $\widehat{{\mathcal M}}$ with $X_+=v \cup \{f,h\}$ and $X_-=\emptyset$. A node $v$ of $K_f$ is a node of $K_{[f,-h]}$ if there is a vector $X$ of $\widehat{{\mathcal M}}$ with $X_+ =  v \cup \{h\}$ and $X_-= f$.  

\begin{lemma}
The component of $K_{[-h,-f]}$ containing the sink of $K_h$ is a monotone path from the sink of $K_h$ to the sink of $K_f$.  The component of $K_{[f,-h]}$ containing the sink of $K_h$ is a monotone path from the source of $K_f$ to the sink of $K_h$.  The only vertex that the paths have in common is the sink of $K_h$.  
\end{lemma}

\begin{proof}

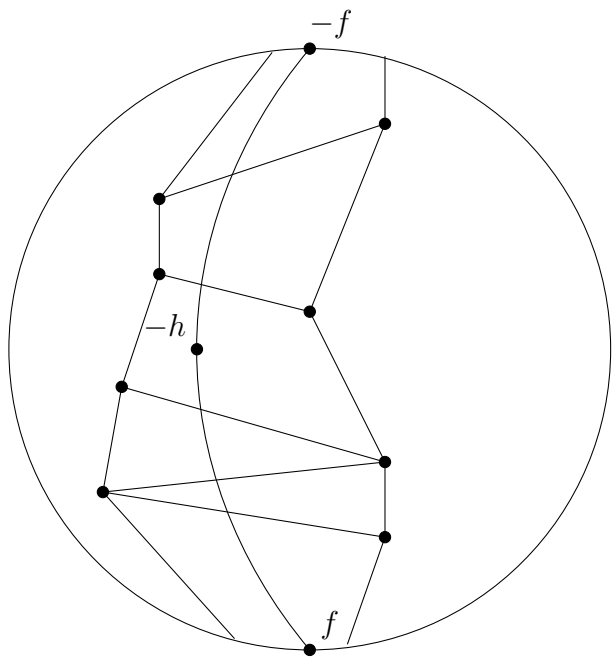
\begin{figure}    
\centering % used for centering Figure
\begin{tikzpicture}[scale=0.5]
\draw circle [radius=8];
    \draw[fill] (-3,0) circle [radius=.15];
    \draw[fill](0,-8) circle [radius=.15];
    \draw[fill](0,8) circle [radius=.15];
    \node[above left] at (-3,0){$-h$};
    \node[above right] at (-0.3,8) {$-f$};
    \node[above right] at (0,-8) {$f$};
%   \draw (0,-8) to [out=-210,in=-150] (0,8);
   \draw (0,-8) to [out=-230,in=-130] (0,8);
 %  \draw (0,-8) to [out=-240,in=-120] (0,8);
    \draw(2,-3)--(-5,-1);
    \draw(-4,2)--(-5,-1);
    \draw(-4,2)--(0,1);
    \draw(-4,2)--(-4,4);
    \draw(2,6)--(-4,4);
    \draw(2,-3)--(-5.5,-3.8);
    \draw(-5.5,-3.8)--(-5,-1);
    \draw[fill] (-4,2) circle [radius=.15];
    \draw[fill] (-5,-1) circle [radius=.15];
    \draw[fill] (-5.5,-3.8) circle [radius=.15];
    \draw(0,1)--(2,6);
    \draw(0,1)--(2,-3);
 %   \draw(2,-3)--(4,-3.8);
    \draw(2,-5)--(2,-3);
    \draw(2,-5)--(-5.5,-3.8);
    \draw(2,-5)--(1,-7.85);
    \draw(-2,-7.7)--(-5.5,-3.8);
    \draw(-1,7.9)--(-4,4);
    \draw(2,6)--(2,7.8);
    \draw[fill] (2,-5) circle [radius=.15];
    \draw[fill] (0,1) circle [radius=.15];
    \draw[fill] (-4,4) circle [radius=.15];
    \draw[fill] (2,6) circle [radius=.15];
    \draw[fill] (2,-3) circle [radius=.15];
\end{tikzpicture}

\caption{Path from source to sink} % title of the Figure
\label{fig:path} % label to refer figure in text
\end{figure}

The digraph $K_{[f,-h]}$ is defined similarly to $K_{[f,h]}$, except that the sign of $h$ is reversed.  This makes the resulting monotone path go to the sink of $K_h$, which is a vertex $v$ where there is a positive vector of $\widehat{{\mathcal M}}$ contained in $v \cup h$.  
In $K_{[-h,-f]}$, the negative signs make the monotone path connect the sinks of $K_h$ and $K_f$, and because $h$ precedes $f$ in $[-h,-f]$, it goes from the sink of $K_h$ to the sink of $K_f$.  If $v$ and $v'$ are adjacent vertices of $K_{[-h,-f]}$, the arc is oriented from $v$ to $v'$ if the cocircuit $Y''$ of $\widehat{{\mathcal M}}$ with zero set $v \cap v'$ has $v\backslash v'$ disagreeing with $h$.  Because $h$ and $f$ have different signs in $Y''$, this means that it has $v\backslash v'$ agreeing with $f$, so the orientation agrees with that of $K_f$.  

 Suppose that a vertex $x$ of $K_f$ is in both $K_{[f,-h]}$ and $K_{[-h,-f]}$. Then there exist vectors $X$ and $X'$ such that $x \cup \{h\} \subseteq X_+$, $ \underline{X}=\underline{X'}=x\cup\{f,h\}$,  $X_f=-,$ $X'_f=+$.  Eliminating $f$ between $X$ and $X'$ shows that $x$ is the sink of $K_h$.
\end{proof}

Concatenation of the component of $K_{[f,-h]}$ containing the sink of $K_h$ and the component of $K_{[-h,-f]}$ containing the sink of $K_h$ yields a monotone path in $K_f$ from the source to the sink.  

Geometric motivation for the path is provided in the realizable case by Figure 3.  
Assume that $\widehat{\mathcal M}$ on $[n] \cup \{f,h\}$ is realized by the $(m+2) \times (n+2)$ matrix 
$$
\widehat{T} = \left[
\begin{array}{c|c|c}
A & 0& 0 \\ \hline
 -c^T& 1&0 \\
-c'^T& 0 & 1
\end{array}\right].
$$
where the last two columns correspond to elements $f$ and $h$. 
Let $\widehat{S}$ be an $(n-m) \times (n+2)$ matrix for which the rows span the orthogonal complement of the row space of $\widehat{T}$.  For Figure 3 we assume that $n-m = r-1 = 3$.  The columns of $\widehat{S}$, scaled  to 
have length 1, are points on the unit sphere.  For every nonnegative cocircuit $Y$ of ${\mathcal M}\backslash f$, the elements of $[n]$ not in the support of $Y$ generate 
a pointed convex cone in $\mathbb{R}^{n-m}$. The cones intersecting the great circle containing $f$ and $h$ correspond to the vertices of the path.

\begin{proposition}
Suppose that ${\mathcal M}$ is an oriented matroid program on $E=[n]\cup\{f,g\}$ and $e\in [n]$ is contained in both the source and the sink of the oriented matroid program.  Then there is a monotone path in $K_f$ from the source to the sink for which no intermediate vertex contains $e$.
\end{proposition}

\begin{proof}
One could let the extension $h$ of the preceding construction be a perturbation of an element $e \in [n]$.  This means that every cocircuit of $\widehat{\mathcal M}$ containing $h$ and $e$ in its support contains both with the same sign.  The property that each of the vectors $X$ encountered in the path, other than the source or sink, has $X_h=+$ implies that 
$X_e\neq +$, so $e \notin v$ for any internal node $v$ of the path.  
\end{proof}

This implies, in particular, that the first digraph of Figure 1 is not the
digraph of an oriented matroid program.  There is a facet $e$, containing the source and the sink of the digraph, for which no monotone path from source to sink has all of its internal nodes off the facet $e$.

\section{Conditions for avoiding a set of circuits}
We want to determine conditions under which the digraph of ${\mathcal M}$ has $r-1$ independent monotone paths from the source to the sink.  We first recall Menger's Theorem:

\begin{proposition} The digraph $K_f$ has $r-1$ independent monotone paths from the source to the sink if and only if there do not exist nodes $\tilde{v}^1, \ldots, \tilde{v}^{r-2}$ of $K_f$, none of which is the source or the sink of $K_f$, such that every monotone path from the source to the sink contains a node of 
$\{\tilde{v}^1, \ldots, \tilde{v}^{r-2}\}$.
\end{proposition}

We will assume that we have a set of nodes $\{\tilde{v}^1, \ldots, \tilde{v}^{r-2}\}$ of $K_f$, none of which is the source or the sink of $K_f$, and determine conditions under which there exists a monotone path from the source to the sink that avoids the nodes in this set.

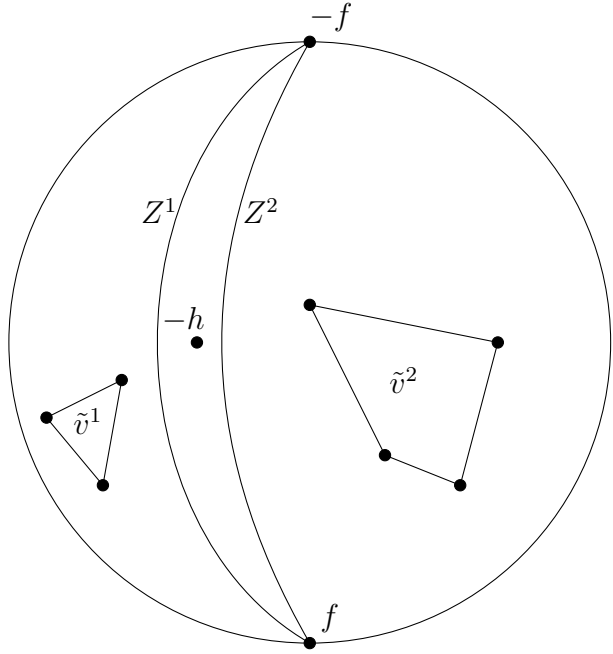
\begin{figure}    
\centering % used for centering Figure
\begin{tikzpicture}[scale=0.5]
\draw circle [radius=8];
    \draw[fill] (-3,0) circle [radius=.15];
    \draw[fill](0,-8) circle [radius=.15];
    \draw[fill](0,8) circle [radius=.15];
    \node[above left] at (-2.5,0){$-h$};
    \node[above right] at (-0.3,8) {$-f$};
    \node[above right] at (0,-8) {$f$};
   \draw (0,-8) to [out=-210,in=-150] (0,8);
    \node at (-4,3.5){$Z^1$};
   \draw (0,-8) to [out=-240,in=-120] (0,8);
    \node at (-1.3,3.5){$Z^2$};
    \draw(-7,-2)--(-5,-1);
    \draw(-7,-2)--(-5.5,-3.8);
    \node at (-5.9,-2.1){$\tilde{v}^1$};
    \draw(-5.5,-3.8)--(-5,-1);
    \draw[fill] (-7,-2) circle [radius=.15];
    \draw[fill] (-5,-1) circle [radius=.15];
    \draw[fill] (-5.5,-3.8) circle [radius=.15];
    \draw(0,1)--(5,0);
    \draw(0,1)--(2,-3);
    \node at (2.5,-1){$\tilde{v}^2$};
    \draw(2,-3)--(4,-3.8);
    \draw(4,-3.8)--(5,0);
    \draw[fill] (0,1) circle [radius=.15];
    \draw[fill] (5,0) circle [radius=.15];
    \draw[fill] (4,-3.8) circle [radius=.15];
    \draw[fill] (2,-3) circle [radius=.15];
\end{tikzpicture}

\caption{Separating $-h$ from cones} % title of the Figure
\label{fig:lnlblock} % label to refer figure in text
\end{figure}

The idea of the argument is illustrated for the realizable case with $r-1=3$ by Fig. 4.    Figure 4 illustrates a typical pair of cones corresponding to facets $\tilde{v}^1,\tilde{v}^2$ of the Las Vergnas sphere of ${\mathcal M}\backslash f$.  The extreme rays of these cones are zero sets of nonnegative cocircuits
$\tilde{Y}^1$ and $\tilde{Y}^2$ of ${\mathcal M}\backslash f$.  These two ($ = r-2$) cones contain neither $f$ nor $-f$, indicating that they do not correspond to the source or the sink of the digraph $K_f$.  Each of the cones corresponding to $\tilde{Y}^1$ and $\tilde{Y}^2$ is separated from point $-h$ by a great circle 
containing points $f$ and $-f$.  These great circles correspond to cocircuits $Z^1$ and $Z^2$ of $\widehat{\mathcal M}$.  Note that there is no point $g$ in the Figure.  The only way that $g$ influences the Figure is in the choices of the subsets of $[n]$ that correspond to supports of nonnegative cocircuits of ${\mathcal M}\backslash f$.

\begin{lemma}\label{avoid} Suppose that $\tilde{v}^1, \ldots, \tilde{v}^{r-2}$ are nodes of $K_f$, but none of them is the sink or source of $K_f$. Suppose that we have an oriented matroid $\widehat{\mathcal M}$ on the set $[n]\cup\{f,h\}$ that has the property that $\widehat{\mathcal M}\backslash h = {\mathcal M}/ g$.
Suppose that for $i=1,\ldots,r-2$ there is a cocircuit $Z^{i}$ of $\widehat{\mathcal M}/ f$ that 
has $Z^i_h=+$ and $Z^i_e \ge 0$ for all $e \in \tilde{v}^i$. Then the monotone path from the source to the sink of $K_f$ created using $\widehat{\mathcal M}$ as in the last section cannot contain any of the vertices $\tilde{v}^1,\ldots,\tilde{v}^{r-2}.$
\end{lemma}

\begin{proof}  If node $\tilde{v}^i$ were on the monotone path, there would be a vector $X$ of $\widehat{{\mathcal M}}$ with $\underline{X} = \tilde{v}^i \cup \{f,h\}$ and $X_e=+$ for $e \in \tilde{v}^i \cup h$. Removing the entry of $X$ indexed by $f$ yields a vector $X'$ of $\widehat{\mathcal M}/f$.  Existence of both the circuit $X'$ and the cocircuit $Z^i$ in $\widehat{\mathcal M}/f$ would violate orthogonality.  
\end{proof}  

\begin{lemma}
Suppose that $\widetilde{\mathcal M}$ is an oriented matroid on $[n] \cup \{h\}$ that is an extension of ${\mathcal M}/\{f,g\}$ by an element $h$ in general position.  Then there exists an extension $\widehat{\mathcal M}$ of ${\mathcal M}/g$ by an element $h$ in general position such that $\widehat{\mathcal M}/f = \widetilde{\mathcal M}$. 
\end{lemma}

\begin{proof}  We define a function $\sigma: {\mathcal C}^*({\mathcal M}/g) \rightarrow \{+,-,0\}$,  and for each $Y \in {\mathcal C}^*({\mathcal M}/g)$ we define a sign vector $\hat{Y}$ on $[n] \cup \{f,h\}$ such that $\hat Y_e = Y_e$ for $e \neq h$ and $\hat Y_h=\sigma(Y)$.  
By Theorem 7.1.8 of \cite{OMbook}, the $\hat Y$ thus defined are cocircuits of an extension of ${\mathcal M}/g$ if they define a single element extension on every rank 2 contraction of ${\mathcal M}/g$.  
For each cocircuit $Y$ of ${\mathcal M}/g$ with $Y_f \neq 0$, define $\sigma(Y) = Y_f$.  For each cocircuit $Y$ of ${\mathcal M}/g$ with $Y_f=0$, the sign vector $Y'$, which is $Y$ with entry $Y_f$ removed, is a cocircuit of ${\mathcal M}/\{f,g\}$.  Let $\hat{Y}'$ be the corresponding cocircuit of $\widetilde{\mathcal M}$ which is the same as $Y'$ except for an extra entry $\hat{Y}'_h$.  Then let $\sigma(Y)=\hat{Y}'_h$. 

This construction is similar to a lexicographic extension.  The rank 2 minors of ${\mathcal M}/g$ contained in ${\mathcal M}/\{f,g\}$ by assumption yield rank 2 minors of $\widehat{\mathcal M}$, while the rank 2 minors of 
${\mathcal M}/g$ not contained in ${\mathcal M}/\{f,g\}$ become rank 2 minors of $\widehat{\mathcal M}$ because $\hat{Y}_h=\hat{Y}_f$ for each cocircuit $Y$ with $Y_f \neq 0$.   
\end{proof}

\begin{corollary}\label{avoid2} Suppose that $\tilde{v}^1, \ldots, \tilde{v}^{r-2}$ are nodes of $K_f$, but none of them is the sink or source of $K_f$. Suppose that we have an oriented matroid $\widetilde{\mathcal M}$ on the set $[n]\cup\{h\}$ that has the property that $\widetilde{\mathcal M}\backslash h = {\mathcal M}/\{f,g\}$.
Suppose that for $i=1,\ldots,r-2$ there is a cocircuit $Z^{i}$ of $\widetilde{\mathcal M}$ that 
has $Z^i_h=+$ and $Z^i_e \ge 0$ for all $e \in \tilde{v}^i$. Then the monotone path from the source to the sink of $K_f$ created using $\widehat{\mathcal M}$ as in the last section cannot contain any of the vertices $\tilde{v}^1,\ldots,\tilde{v}^{r-2}.$
\end{corollary}

\begin{lemma}\label{r-2} For $i=1,\ldots,r-2$ and each $e \in \tilde{v}^i$ there exists a cocircuit $Y^{e,i}$ of ${\mathcal  M}/ \{f,g\}$ such that $Y^{e,i}_e = +$, $Y^{e,i}$ is nonnegative on $\tilde{v}^i$.  Composition of the $Y^{e,i}$ for a given $i$ yields a covector $Y^i$ of ${\mathcal M}/\{g,f\}$ that is positive on $\tilde{v}^i$.
\end{lemma}
\begin{proof} If there were no such cocircuit $Y^{e,i}$, then there would be a circuit $X^i$ of ${\mathcal M}/g$ with support in $\tilde{v}^i \cup \{f\}$ that is nonnegative on $\tilde{v}^i$. (See Corollary 10.2.18 of \cite{OMbook}.)  Corresponding to $X^i$ there is a vector $\overline{X}^i$ of ${\mathcal M}$ such that $X^i$ is obtained from $\overline{X}^i$ by removing entry $g$.  The support of $\overline{X}^i$ is not contained in $[n]$, by acyclicity of ${\mathcal M}$. There is a covector $Y^i$ of ${\mathcal M}$ for which $\tilde{v}^i$ is the zero set.  The supports of $\overline{X}^i$ and $Y^i$ intersect, if at all, in a subset of $\{f,g\}$.  It is assumed that $g \in \underline{Y}^i$,  so the support of $\overline{X}^i$ must contain $f$ as well as $g$.  Therefore, $X^i_f \neq 0$.  This is not possible due to the assumption that $\tilde{v}^i$ is neither the source nor the sink of $K_f$.
\end{proof}

We now show our main result under the assumption that ${\mathcal M}/\{f,g\}$ has an {\it adjoint}.  
 Existence of an adjoint means that the oriented matroid can be faithfully depicted by a drawing such as those in Figures 3 and 4.  Each hyperplane of the oriented matroid can be replaced by a pseudohyperplane such that the resulting set of regions of $\mathbb{R}^{r-2}$ are the covectors of the adjoint oriented matroid.

\begin{proposition}
Suppose that $\tilde{v}^1, \ldots, \tilde{v}^{r-2}$ are nodes of $K_f$, but none of them is the sink or source of $K_f$. Suppose that the oriented matroid ${\mathcal M}/\{f,g\}$ has an adjoint.  Then there exists an extension $\widetilde{\mathcal M}$ of ${\mathcal M}/\{f,g\}$ by an element $h$ such that for $i=1,\ldots,r-2$ there is a cocircuit $Z^{i}$ of $\widetilde{\mathcal M}$ that 
has $Z^i_h=+$ and $Z^i_e \ge 0$ for all $e \in \tilde{v}^i$.
\end{proposition}

\begin{proof}
Suppose that ${\mathcal M}/\{f,g\}$ has an adjoint ${\mathcal A}$. By Definition 5.3.5 of \cite{OMbook}, ${\mathcal A}$ is an oriented matroid, also of rank $r-2$, and there is an injection $\alpha$ from the ground set $E^{ad}$ of ${\mathcal A}$ to the set of cocircuits of ${\mathcal M}$, so that the image of $\alpha$ contains exactly one element of each pair of opposite cocircuits of ${\mathcal M}/\{f,g\}$.  For each element $e$ of the ground set of ${\mathcal M}/\{f,g\}$, there is a cocircuit $Z^e$ of ${\mathcal A}$ with $Z^e_{e'}=\alpha(e')_e$ for all $e' \in E^{ad}$.

Each of the sets $\tilde{v}^1, \ldots, \tilde{v}^{r-2}$ is a set of rank $r-2$ in ${\mathcal M}/\{f,g\}$. For $i=1,2,\ldots,r-2$, let $E^i$ be the set of elements $e'$ of $E^{ad}$ for which $\alpha(e')$ is nonnegative on $\tilde{v}^i$ or nonpositive on $\tilde{v}^i$. Because $\tilde{v}^i$ has rank $r-2$ in ${\mathcal M}/\{f,g\}$, it follows from Lemma \ref{r-2} that $E^i$ has rank $r-2$ in ${\mathcal A}$.  Let $B$ be a basis for ${\mathcal A}$ obtained by taking elements $(e')^i \in E^i$, for $i=1,2,\ldots,r-2$.  

Because $B$ is a basis for ${\mathcal A}$, there exists a covector $H$ of ${\mathcal A}$ such that $H_{(e')^i}=+$ for each $i$ for which $\alpha((e')^i)$ is nonnegative on $\tilde{v}^i$ and $H_{(e')^i}=-$ for each $i$ for which $\alpha((e')^i)$ is nonpositive on $\tilde{v}^i$.  By Theorem 7.5.8 of \cite{OMbook}, the 
covector $H$ implies the existence of an extension $\tilde{\mathcal M}$  of ${\mathcal M}/\{f,g\}$ by an element $h$ so that for each $i$ with $\alpha((e')^i)$ nonnegative on $\tilde{v}^i$ we have $\alpha((e')^i)_h=+$ and for each $i$ with $\alpha((e')^i)$ nonpositive on $\tilde{v}^i$ we have $\alpha((e')^i)_h=-.$  For $i$ in the first case, let $Z^i$ be $\alpha((e')^i)$, and in the second case, let $Z^i$ be $-\alpha((e')^i)$.
\end{proof}

This proof was inspired by the proof of Theorem 4.3 in \cite{GMR}.  The proof allows for great flexibility in constructing the basis $B$.  On the other hand, the assumption that ${\mathcal M}/\{f,g\}$ has an adjoint is a very strong assumption in oriented matroid theory.  

\begin{corollary}\label{Zsexist}Suppose that $\tilde{v}^1, \ldots, \tilde{v}^{r-2}$ are nodes of $K_f$, but none of them is the sink or source of $K_f$. Suppose that the oriented matroid ${\mathcal M}/\{f,g\}$ has rank $r-2 \le 3$.  Then there exists an extension $\widetilde{\mathcal M}$ of ${\mathcal M}/\{f,g\}$ by an element $h$ such that for $i=1,\ldots,r-2$ there is a cocircuit $Z^{i}$ of $\widetilde{\mathcal M}$ that 
has $Z^i_h=+$ and $Z^i_e \ge 0$ for all $e \in \tilde{v}^i$.
\end{corollary}

\begin{proof}
Proposition 6.3.6 of \cite{OMbook} states that every oriented matroid of rank 3 has an adjoint.  
\end{proof}

\begin{theorem} The Holt-Klee condition holds for rank 4 and rank 5 oriented matroid programs, generalizing the Holt-Klee condition for 3- and 4-dimensional polytopal digraphs.  
\end{theorem}

\begin{proof}The Theorem follows from Menger's Theorem, Lemma \ref{avoid} and Proposition \ref{Zsexist}
\end{proof}

\section{A weaker condition on ${\mathcal M}/\{f,g\}$}  

In the previous section, we were given nodes $\tilde{v}^1,\tilde{v}^2,\ldots,\tilde{v}^{r-2}$ of $K_f$, none of which are the source and sink of $K_f$, and are asked if there is an extension $\widetilde{\mathcal M}$ of ${\mathcal M}/\{f,g\}$ by an element $h$ such that for $i=1,2,\ldots,r-2$ there is a cocircuit $Z^i$ of $\widetilde{\mathcal M}$ such that $Z^i_h=+$ and $Z^i_e \ge 0$ for all $e \in \tilde{v}^i$.  

We would like to compare the existence of such an extension to {\it Levi's Intersection Property} (See Definition 7.5.2 of \cite{OMbook}):  Given an oriented matroid ${\mathcal N}$ of rank $r-2$ and cocircuits 
$Y^1,Y^2,\ldots,Y^{r-3}$ of ${\mathcal N}$, there exists an extension $\widehat{\mathcal N}$ of ${\mathcal N}$  by an element $h$ and 
cocircuits $Z^1,Z^2,\ldots,Z^{r-3}$ of $\widehat{\mathcal N}$ such that for $i=1,2,\ldots,r-3$, $Z^i_h=0$ and $Z^i_e =Y^i_e$ for all $e \neq h$.  

The two conditions differ in three major ways: 
\begin{itemize}
    \item Levi's intersection property requires $r-3$ cocircuits $Z^i$ rather than the $r-2$ of the previous section.
    \item The cocircuits $Y^1,Y^2,\ldots,Y^{r-3}$ in Levi's intersection property are replaced by rank $r-2$ sets $\tilde{v}^1, \tilde{v}^2,\ldots,\tilde{v}^{r-2}$.  
    \item The cocircuits $Z^i$ of the extension of the previous section must have $Z_h=+$, rather than $Z_h=0$ as in Levi's intersection property.    
\end{itemize}
  
We introduce a generalized Levi's intersection property: 
\begin{definition}
Suppose ${\mathcal N}$ is an oriented matroid of rank $r-2$.  ${\mathcal N}$ satisfies the generalized Levi's intersection property if for any set of covectors $\{Y^1,Y^2,\ldots,Y^{r-3}\}$, there exists an extension $\widehat{\mathcal N}$ of ${\mathcal N}$ by an element $h$ and covectors $W^1,W^2,\ldots,W^{r-3}$ of $\widehat{\mathcal N}$ such that for $i=1,2,\ldots,r-3$, $W^i_h=0$ and $W^i_e =Y^i_e$ for all $e \neq h$.   
\end{definition}

\begin{proposition}
If ${\mathcal M}/\{f,g\}$ satisfies the generalized Levi's intersection property, then ${\mathcal M}$ has the Holt-Klee property.
\end{proposition}

\begin{proof}
Suppose $\tilde{v}^1,\tilde{v}^2,\ldots,\tilde{v}^{r-2}$ are nodes of $K_f$ and none is the source or the sink of $K_f$.  For $i=1,2,\ldots,r-2$, Lemma \ref{r-2} produces a covector $Y^i$ of ${\mathcal M}/\{f,g\}$ with $\tilde{v}^i\subseteq Y^i_+$.  By the generalized Levi's intersection property, there exists an extension $\widetilde{\mathcal M}$ of ${\mathcal M}/\{f,g\}$ by an element $h$ and covectors $W^1,W^2,\ldots,W^{r-3}$ of $\widetilde{\mathcal M}$ such that for $i=1,2,\ldots,r-3$, $W^i_h=0$ and $W^i_e =\tilde{Y}^i_e$ for all $e \neq h$. 

Let $\tilde{Y}^{r-2}$ be a covector of $\widetilde{\mathcal M}$ that agrees with $Y^{r-2}$ on $[n]$.  If $\tilde{Y}^{r-2}_h=0$, we can replace $\tilde{Y}^{r-2}$ by the composition of $\tilde{Y}^{r-2}$ with a covector of $\widetilde{\mathcal M}$ that is positive on element $h$.  If $\tilde{Y}^{r-2}_h=-$, get a new extension by reversing element $h$ in $\widetilde{\mathcal M}$.  Thus we can assume that $\tilde{v}^{r-2}\cup \{h\} \subseteq Y^{r-2}_+$.  By Proposition 3.7.2 of \cite{OMbook}, there is a cocircuit $Z^{r-2}$ of $\widetilde{\mathcal M}$ such that $e \in \underline{Z}^{r-2}$ and $Z^{r-2}$ agrees with $Y^{r-2}$ wherever it is nonzero.  Because the rank of $\tilde{v}^{r-2}$ is $r-2,$  $\underline{Z}^{r-2}\cap \tilde{v}^{r-2} \neq \emptyset.$  

Note that even if $h$ was negated in the construction of $Z^{r-2}$, the vectors $W^1,W^2,\ldots W^{r-3}$ are still covectors of $\widetilde{\mathcal M}$.  For each $i =1,2,\ldots,r-3$, We can replace each $W^i$ with the composition of $W^i$ with a covector of $\widetilde{\mathcal M}$ that is positive on element $h$.  As before, there is a cocircuit $Z^i$ of $\widetilde{\mathcal M}$ that is positive on element $h$ and agrees with $W^i$ on $\underline{Z}^i$.  

The existence of the extension implies that ${\mathcal M}$ has the Holt-Klee property, as in the last section.

\end{proof}

\section{Non-lifting subdivisions and P-cubes} 

The digraph $K_f$ of an oriented matroid program can be derived from the pair $({\mathcal M}/g, {\mathcal V})$, where ${\mathcal V}$ is the set of vertices of the Edmonds-Mandel sphere of ${\mathcal M}\backslash f$.  In this section, we work with a similar pair $({\mathcal N} \cup f,S)$, where $S$ is a collection of subsets of the ground set of a rank $r-1$ oriented matroid ${\mathcal N}$, but we will not assume that the members of this pair are obtained from a larger oriented matroid ${\mathcal M}$.  The role of the oriented matroid ${\mathcal M}\backslash{f}/g$ will in this section be played by the oriented matroid ${\mathcal N}$. Recall Definition \ref{defsubdivision}, taken from section 4 of \cite{S}.

A subdivision $S$ of an oriented matroid $N$ is called a {\it lifting} subdivision if there is an oriented matroid $\widehat{\mathcal N}$ on a set $E \cup g$ so that ${\mathcal N}= \widehat{\mathcal N}/g$ and the subsets in $S$ are the complements of the supports of nonnegative cocircuits of $\widehat{\mathcal N}$.  

We will also make the following assumption on ${\mathcal N}$: There are no nonnegative cocircuits (facets) of ${\mathcal N}$.    This implies in the realizable case that subsets of the columns of a matrix realizing ${\mathcal N}$ indexed by elements of $S$ form a complete pointed fan.  

A subdivision defines a graph $G$ for which the vertices are the cells of $S$, and two cells $\sigma_1$ and $\sigma_2$ are adjacent if the rank of $\sigma_1 \cap \sigma_2$ has rank $r-2$.  
An extension of ${\mathcal N}$ by an element $f$ in general position defines an orientation $K_f$ of $G$.  Whenever two cells $\sigma_1$ and $\sigma_2$ are adjacent, there is a cocircuit $Y$ of ${\mathcal N} \cup f$ for which $Y_e=0$ for $e \in \sigma_1 \cap \sigma_2$ and $Y_f=+$.  The edge is directed from $\sigma_1$ to $\sigma_2$ if $Y_e=+$ for all $e \in \sigma_1 \backslash \sigma_2$.  

The paper \cite{S} notes that the definition of subdivision implies: 
Every extension ${\mathcal N}\cup f$ in general position is covered by exactly one cell of $S$.  This means that there is a unique cell $\sigma_1$ such that 
there is a vector $C^1$ for which $C^1_f=-$ and $C^1_+=\sigma_1$, and there is a unique cell $C^2$ for which $C^2_f=+$ and $C^2_+=\sigma_2$.  The cell $\sigma_1$ is the source of the digraph $K_f$ and $\sigma_2$ is the sink.  

The monotone paths constructed in Section 4 are defined in the greater generality of subdivisions in \cite{S}, using Lemma \ref{Gfhdef}

The arguments of Section 5 only use the properties of $K_f$ implied by the definition of a subdivision of an oriented matroid, so we get the following.
\begin{proposition} Let ${\mathcal N}$ be an oriented matroid of rank $r-1$ on a set $E$, and suppose that ${\mathcal N}$ has no nonnegative cocircuits.  
Suppose $K_f$ is a digraph obtained from a subdivision of ${\mathcal N}$ and a one-element extension ${\mathcal N} \cup f$ of ${\mathcal N}$.  If $r \le 5$, then $K_f$ has $r-1$ independent monotone paths from source to sink.   
\end{proposition}

An application of this proposition concerns P-oriented matroid complementarity problems, introduced by Todd \cite{Todd} and studied by Klaus and Miyata \cite{KM}.  Here ${\mathcal N}$ is a rank $n$ oriented matroid on the set $[2n]$, satisfying the condition: 

\begin{equation}\label{eqn:einstein}
\mbox{For each circuit } Y \mbox{ of } {\mathcal N}, \mbox{ there is } i \in \{1,2,\ldots,n\}\mbox{ with }Y_i = Y_{n+i} \neq 0. \tag{P}
\end{equation}

If ${\mathcal N}$ satisfying the condition above is realized by an $n \times 2n$ matrix of the form $[M,I]$, then the matrix $M$ must be a P-matrix, i.e. its principal minors are all positive (see \cite{Todd}). If ${\mathcal N}$ satisfies property (P), then a subdivision of ${\mathcal N}$ is given by $S=\{B \subseteq [2n]: |B \cap \{i,n+i\}|=1 \mbox { for } i=1,2,\ldots,n\}.$ The undirected graph of such a subdivision is that of an $n$-cube.  An extension of the oriented matroid realized by an $n \times 2n$ matrix $[M,I]$ gives rise to a {\it P-matrix linear complementarity problem,} which is to find the unique sink in the digraph $K_f$ defined by the subdivision and the extension.

In \cite{GMR}, it was shown that for a realizable extension of an oriented matroid satisfying property (P), the Holt-Klee property holds for the resulting digraph.

If we require ${\mathcal N}$ to satisfy property $P$, but do not require the extension ${\mathcal N} \cup f$ to be realizable, we get what Klaus and Miyata \cite{KM} call a POMCP, where OMCP stands for oriented matroid complementarity problem.  Our results imply that for such POMCP with the rank of ${\mathcal N}$ at most 4, the digraph satisfies the Holt-Klee property.  Klaus and Miyata were able to enumerate all the 6910 combinatorially distinct orientations $K_f$ of the $4$-cube obtained from POMCPs.  They found the surprising fact that each of these 6910 orientations was obtainable from a {\it realizable} oriented matroid.   

The example of Section 4.4 of $\cite{GMR}$ can be shown by the same argument given there to be an example of an orientation of the 4-cube that has a unique source and sink on every face and satisfies the Holt-Klee condition but is not given by an POMCP.   As POMCP orientations of the cube are, in theory, more general than those obtained from oriented matroid programs with Edmonds-Mandel spheres combinatorially equivalent to the cube, this example also shows that the Holt-Klee property, together with the condition that there is a unique source and unique sink for each face, are not sufficient to imply that an orientation is obtainable from an oriented matroid program.  

The existing literature does not yet contain an example of a POMCP for which the digraph is the digraph of an oriented matroid program.  There also appears to be no example known of an acyclic realizable POMCP digraph that cannot be obtained from a linear program over a cube.

\section{Rank 6 oriented matroid programs} If ${\mathcal M}$ is a rank 6 oriented matroid on $[n] \cup \{f,g\}$, and we are given 4 feasible cocircuits, none of which is the source or the sink of $K_f$, then we do not know if there exists a path from the source of $K_f$ to the sink that avoids all four given cocircuits.  It seems doubtful that an extension $h$ of ${\mathcal M}/\{f,g\}$ guaranteeing such a path, as in section 5, can always be found.  However, there may be completely different approaches to proving that the Holt-Klee property holds, if it indeed does. As inspiration, one could consider the proof of Balinski's theorem due to Barnette \cite{Bar}.

%Suppose, for example, that the minor ${\mathcal M}/\{f,g\}$ is the rank 4 oriented matroid $RS(8)$ on 8 elements.  It is known that for $RS(8)$ there is an element $e*$ with the property that there are only three extreme cocircuits of $RS(8)$ on the positive side of 
%element $d^*$.  
%For each of these three extreme cocircuits$Y^i$ for $i=1,2,3$, let $T^i$ be the simplicial region between $Y^i$ and element $e^*$.  
%Suppose $T^i$ is the tope that is negative on the elements of $[n]\backslash \underline{\tilde{Y^i}}$ for $i = 1,2,3$.  If the pseudosphere 
%corresponding to the new element $h$ intersected the interior of each of the topes $T^i$, then it follows that all of the cocircuits that
%are positive on component $e*$ have the same sign, say $+$, on element $e*$,  IF the tope that is negative on the elements of 
%$[n]\backslash \underline{\tilde{Y^4}}$ were a tope on the negative side of $e*$, then it would not be possible for a vertex of this tope to be positive on element $h$.  

The author would like to express his appreciation for the careful reading and helpful comments of the referee.
A preliminary version of this material was first presented in a talk at Casa Matematic\'as Oaxaca in November, 2015.

http://videos.birs.ca/2015/15w5006/201511021203-MorrisJr..mp4

\end{document}